\documentclass[reqno,10pt]{amsart}

\usepackage{amssymb}
\usepackage[all,cmtip]{xy}
\newtheorem{Proposition}{Proposition}[section]
\newtheorem{MainTheorem}{Theorem}
\newtheorem{Definition}[Proposition]{Definition}

\newtheorem{Theorem}[Proposition]{Theorem}
\newtheorem{Corollary}[Proposition]{Corollary}

\DeclareMathOperator{\vol}{vol}
\DeclareMathOperator{\vp}{vp}
\newcommand{\R}{\mathbb{R}}

\title{The isoperimetrix in the dual Brunn-Minkowski theory}

\author{Andreas Bernig}
\email{bernig@math.uni-frankfurt.de}
\address{Institut f\"ur Mathematik, Goethe-Universit\"at Frankfurt, Robert-Mayer-Str. 10, 60054
Frankfurt, Germany}

\thanks{ Supported by DFG grants BE 2484/3-1 and BE 2484/5-1.\\ AMS 2010 {\it Mathematics subject classification}:
52A38, 
52A40, 
53C60. 
}
\begin{document}

\begin{abstract}
We introduce the dual isoperimetrix which solves the
isoperimetric problem in the dual Brunn-Minkowski theory. We then show how the dual isoperimetrix is related to the isoperimetrix from the Brunn-Minkowski theory.  
\end{abstract}

\maketitle

\section{Introduction and statement of main results} 

A {\it definition of volume} is a way of measuring volumes on finite-dimensional normed spaces, or, more
generally, on Finsler manifolds. Roughly speaking, a definition of
volume $\mu$ associates to each finite-dimensional normed space $(V,\|\cdot\|)$ a norm on $\Lambda^n V$, where $n=\dim
V$. We refer to \cite{alvarez_thompson, thompson_book96} and Section \ref{sec_defvol} for more information. The
best known examples of definitions of volume are the {\it Busemann volume} $\mu^b$, which equals the Hausdorff measure,
the {\it Holmes-Thompson volume} $\mu^{ht}$, which is related to symplectic geometry and {\it Gromov's mass*} $\mu^{m*}$, which thanks to its convexity properties is often used in geometric measure theory. 

To each definition of volume $\mu$ may be associated a dual definition of volume $\mu^*$. For instance, the dual of
Busemann's volume is Holmes-Thompson volume and vice versa. The dual of Gromov's mass* is Gromov's mass, which however
lacks good convexity properties and is used less often. 

Given a definition of volume $\mu$ and a finite-dimensional normed vector space $V$, there is an induced
$(n-1)$-density $\mu_{n-1}:\Lambda^{n-1} V \to \R$. Such a density may be integrated over $(n-1)$-dimensional
submanifolds in $V$.

Given a compact convex set $K \subset V$, we let 
\begin{displaymath}
 A_\mu(K):=\int_{\partial K} \mu_{n-1}
\end{displaymath}
be the $(n-1)$-dimensional surface area of $K$. 

If the boundary $\partial K$ is smooth, each tangent space $T_p \partial K \subset V$ carries the induced norm and
$\partial K$ is a Finsler manifold. In the general case, one may make sense of the integral by using Alexandrov's
theorem \cite{alexandrov39}.  

The definition of volume $\mu$ is called {\it convex}, if for compact convex bodies $K \subset L$, we have $A_\mu(K)
\leq
A_\mu(L)$.
There are many equivalent ways of defining convexity of volume definitions, we refer to \cite{alvarez_thompson} for
details. The above mentioned three examples are convex. 

Given a convex definition of volume and an $n$-dimensional normed space $V$ with unit ball $B$, there is a unique
centrally symmetric compact convex body $\mathbb{I}_\mu B$ such that 
\begin{displaymath}
 A_\mu(K)=n V(K[n-1],\mathbb{I}_\mu B), \quad K \in \mathcal{K}(V).
\end{displaymath}
Here $V$ denotes the mixed volume and $\mathcal{K}(V)$ stands for the space of compact convex bodies in $V$.
$\mathbb{I}_\mu B$ is called the isoperimetrix \cite{alvarez_thompson}. It was introduced in 1949 by Busemann
\cite{busemann49b} and has applications in crystallography \cite{taylor78, wulff01} and in
geometric measure theory \cite{ambrosio_kirchheim00, morgan92, taylor73}.

As its name indicates, the isoperimetrix is related to isoperimetric problems. More precisely, Busemann
\cite{busemann49b} showed that among all compact convex bodies of a given, fixed volume, a homothet of the
isoperimetrix has minimal surface area. 

The isoperimetrices of the above mentioned examples of definitions of volumes are related to important concepts from
convex geometry. For Busemann's definition of volume, we have $\mathbb{I}_{\mu^b} B=\omega_{n-1} (IB)^\circ$, where
$IB \subset V^*$ is the {\it intersection body} of $B$, $I^\circ B:=(IB)^\circ \subset V$ its polar body, and
$\omega_{n-1}$ the volume
of the $(n-1)$-dimensional (Euclidean) unit ball. For the Holmes-Thompson volume, we have
$\mathbb{I}_{\mu^{ht}}B=\frac{1}{\omega_{n-1}}
\Pi(B^\circ)$, where $\Pi$ denotes the {\it projection body}. It was shown recently by Ludwig \cite{ludwig10} that the
Holmes-Thompson surface area can be uniquely characterized by a valuation property. The
isoperimetrix for Gromov's mass* is a dilate of the
{\it wedge body} of $B$. 

In the present paper, we introduce a {\it dual isoperimetrix} which belongs to the dual Brunn-Minkowski theory.
In the dual Brunn-Minkowski theory, the natural setting is that of star bodies (i.e. compact, star shaped bodies
containing
the origin in their interior, with continuous radial function). The Minkowski sum of
convex bodies is replaced by the radial sum and mixed volumes by dual mixed volumes. The dual Brunn-Minkowski theory
was developed by Lutwak \cite{lutwak75, lutwak88} and plays a prominent role in modern convexity.

Let us describe our main results. The space of star bodies with smooth radial function in a finite-dimensional vector
space $V$ is denoted by
$\mathcal{S}^\infty(V)$. Let $S \in \mathcal{S}^\infty(V)$ and $p \in \partial S$. We will write $\langle
p\rangle:=\R \cdot p$ for the line generated by $p$. The tangent space $T_p \partial
S$ is naturally identified with the quotient space $V/\langle
p\rangle$ and inherits the quotient norm. Therefore $\partial S$ is in a natural way a Finsler manifold. Note that
this metric is not the induced metric in general. It was (to the best of our knowledge) first studied in a recent
paper by Faifman \cite{faifman}. 

Let $\mu$ be a definition of volume, $V$ an $n$-dimensional normed vector space and $M$ a smooth $(n-1)$-dimensional
submanifold such that $p \notin T_pM$ for all $p \in M$. Then $T_pM \cong V/\langle p\rangle$ can be endowed with the quotient norm and we obtain an $(n-1)$-density $\tilde \mu_{n-1}:\Lambda^{n-1} TM \to
\R$. In the particular case where $M=\partial S$ for some smooth star-shaped set $S$, we denote by 
\begin{displaymath}
 \tilde A_\mu(S):=\int_{\partial S} \tilde \mu_{n-1}
\end{displaymath}
the surface area of $S$ with respect to $\mu$ and the quotient metric (see Section \ref{sec_dual_isoperimetrix} for more
details). By continuity in the radial topology, this extends to non smooth star bodies as well.

Let us now formulate our main theorem. 

\begin{MainTheorem} \label{mainthm_dual_isop}
Let $\mu$ be a definition of volume, and let $V$ be an $n$-dimensional normed space with unit ball $B$. 
\begin{enumerate}
 \item There
exists a star body $\tilde{\mathbb{I}}_\mu B \subset V$ such that 
\begin{displaymath}
 \tilde A_\mu(S)=n \tilde V(S[n-1],\tilde{\mathbb{I}}_\mu B), \quad S \in \mathcal{S}^\infty(V). 
\end{displaymath}
Here $\tilde V$ denotes the dual mixed volume \cite{lutwak75}. We call  $\tilde{\mathbb{I}}_\mu B$ the \emph{dual isoperimetrix}. 
\item Dual isoperimetric problem: Among all star bodies of a given volume, a
dilate of the dual isoperimetrix has \emph{maximal} surface area.
\item Suppose that the dual definition of volume $\mu^*$ is convex. Then the usual
isoperimetrix for the dual
definition of volume $\mu^*$ and the dual isoperimetrix are related by 
\begin{displaymath}
 \tilde{\mathbb{I}}_\mu B=\left(\mathbb{I}_{\mu^*} B^\circ\right)^\circ.
\end{displaymath}
\end{enumerate}
\end{MainTheorem}

\begin{Corollary}
The dual isoperimetrix for Busemann's definition of volume $\mu^b$ is  
\begin{displaymath}
 \tilde{\mathbb{I}}_{\mu^b}(B)=\omega_{n-1} \Pi^\circ(B),
\end{displaymath}
while for the Holmes-Thompson volume $\mu^{ht}$, we have 
\begin{displaymath}
 \tilde{\mathbb{I}}_{\mu^{ht}}(B)=\frac{1}{\omega_{n-1}} I(B^\circ).
\end{displaymath}
\end{Corollary}

Our second main theorem is an affinely invariant inequality. A survey over affinely invariant inequalities can be found in
\cite{lutwak_handbook93}. 

\begin{MainTheorem}[Surface area of the unit sphere] \label{mainthm_surface_sphere}
Let $(V,B)$ be an
$n$-dimensional normed space. Then 
\begin{equation} \label{eq_dual_affine_ineq}
\tilde A_{{\mu}^b}(B) \leq
n \omega_n,
\end{equation} 
Equality is attained precisely for centered ellipsoids.
\end{MainTheorem}

In the two-dimensional case, this bound was conjectured by Faifman \cite{faifman}, who gave the
non-optimal upper bound of $8$. He also gave a lower bound of $4$ and conjectured that $8 \log 2 \approx 5.5$ (which appears in the case of a square)  
is the optimal lower bound. Using John's ellipsoid, one can improve the lower bound in the two-dimensional case to $\sqrt{2} \pi \approx 4.4$, but we do not have a sharp lower bound.  

As a corollary, we prove an upper bound for the quotient girth. Recall that the {\it girth} of a normed space is the
length of the shortest symmetric curve on the unit sphere, measured with the Finsler metric induced by the norm.
Analogously, the quotient girth is the length of the shortest symmetric curve on the unit sphere, measured with the
quotient Finsler metric. 

\begin{Corollary} \label{cor_max_quot_girth}
In any dimension, the quotient girth is bounded from above by $2\pi$, with equality precisely for ellipsoids. 
\end{Corollary}

\subsubsection*{Acknowledgements}
I wish to thank Dmitry Faifman, Monika Ludwig, Franz Schuster and Deane Yang for interesting discussions and useful
remarks on a
first draft of
this paper.
 
\section{The isoperimetrix in the Brunn-Minkowski theory}
\label{sec_defvol}

In this section, we recall some classical notions from convex geometry and Finsler geometry. References for this
section are \cite{alvarez_thompson, gardner_book06, thompson_book96}. 

\begin{Definition}[Definition of volume]
A definition of volume assigns to each normed vector space $(V,\|\cdot\|)$ of dimension $n<\infty$ a norm $\mu_V$ on
the one-dimensional space $\Lambda^nV$ such that the following three conditions are satisfied:
\begin{enumerate}
 \item If $V$ is Euclidean, then $\mu_V$ is the norm induced by usual Lebesgue measure.
\item If $f:(V,\|\cdot\|) \to (W,\|\cdot\|)$ is a linear map that does not increase distances, then the induced map
$\Lambda^nf:(\Lambda^nV,\mu_V) \to (\Lambda^nW,\mu_W)$ does not increase distances. 
\item The map $(V,\|\cdot\|) \mapsto (\Lambda^n V,\mu_V)$ is continuous with respect to Banach-Mazur distance. 
\end{enumerate}
If $B$ is the unit ball of $V$, we also write $(V,B)$ instead of $(V,\|\cdot\|)$.
\end{Definition}

We will often need an alternative description of a definition of volume. Let $\mu$ be a definition of volume and
$(V,B)$ a normed space. Set $\mathcal{V}(B):=\int_B \mu$. Then $\mathcal{V}$ is a functional on the space
of compact convex, centrally symmetric bodies with the origin in their
interior. It is continuous, invariant under invertible linear maps and satisfies the following monotonicity condition:
whenever $B_1 \subset B_2$ belong to the same $n$-dimensional vector space, then 
\begin{displaymath}
 \frac{\mathcal{V}(B_1)}{\mathcal{L}^n B_1} \geq \frac{\mathcal{V}(B_2)}{\mathcal{L}^n B_2}.
\end{displaymath}
Here $\mathcal{L}^n$ is any choice of Lebesgue measure. Conversely, any continuous invariant functional with this
monotonicity property comes from a definition of volume.

Before presenting some examples, we have to recall some notions from convex geometry. Given a compact convex body $K
\subset V$, the support function is defined by 
\begin{displaymath}
 h_K(\xi):=\sup_{x \in K} \xi(x), \quad \xi \in V^*. 
\end{displaymath}
This function is positively $1$-homogeneous and convex. Conversely, any positively $1$-homogeneous and convex function
is the support function of a unique compact convex body. 

The polar body is defined by
\begin{displaymath}
 K^\circ:=\{\xi \in V^*: \langle \xi,x\rangle \leq 1, \quad \forall x \in K\}.
\end{displaymath}

The radial function of a compact convex set $K$ containing the origin in its interior is given by
\begin{displaymath}
 \rho_K(x):=\max\{\lambda \geq 0: \lambda x \in K\}, \quad x \in V \setminus \{0\}. 
\end{displaymath}
We will make frequent use of the well-known identity
\begin{equation} \label{eq_relation_support_radial}
 \rho_K(x)=\frac{1}{h_{K^\circ}(x)}.
\end{equation}

Examples 
\begin{enumerate}
 \item The {\it Busemann volume} is the definition of volume such that $\mathcal{V}(B)=\omega_n$ for any $n$-dimensional
unit ball $B$. Here $\omega_n$ is the Euclidean volume of the Euclidean unit ball.   
\item The {\it Holmes-Thompson volume} is defined by 
\begin{displaymath}
 \mathcal{V}(B)=\frac{1}{\omega_n} \vp(B),
\end{displaymath}
where $\vp(B)$ denotes the volume product of $B$, i.e. the symplectic volume of $B \times B^\circ \subset V \times
V^*$. 
\item \label{item_gromov} Gromov's mass is defined by 
\begin{displaymath}
\mathcal{V}(B):= \frac{2^n}{n!} \frac{\mathcal{L}^n(B)}{\sup_{C \subset B} \mathcal{L}^n(C)}. 
\end{displaymath}
Here $\mathcal{L}^n$ is any Lebesgue measure and $C$ ranges over all cross-polytopes inscribed in $B$. 
\item \label{item_benson} Benson definition (or Gromov's mass*). 
\begin{displaymath}
\mathcal{V}(B):= 2^n \frac{\mathcal{L}^n(B)}{\inf_{P \supset B} \mathcal{L}^n(P)}. 
\end{displaymath}
Here $P$ ranges over all parallelotopes circumscribed to $B$.
\item  \label{item_ivanov} Ivanov's definition of volume \cite{ivanov08}
\begin{displaymath}
\mathcal{V}(B):= \omega_n \frac{\mathcal{L}^n(B)}{\mathcal{L}^n(E)},
\end{displaymath}
where $E$ is the maximal volume ellipsoid inscribed in $B$ (i.e. the John ellipsoid). 
\item \label{item_dual_ivanov} The dual of Ivanov's definition of volume
\begin{displaymath}
\mathcal{V}(B):= \omega_n \frac{\mathcal{L}^n(B)}{\mathcal{L}^n(E)}, 
\end{displaymath}
where $E$ is the minimal volume ellipsoid circumscribed to $B$ (i.e. the L\"owner ellipsoid). 
\end{enumerate}
 
If $\mu$ is a definition of volume and $V$ an $n$-dimensional vector space, $\mu_{V^*}$ is a norm on the one-dimensional
vector space $\Lambda^n V^*$, whose dual $(\mu_{V^*})^*$ is a norm on $\Lambda^n V$. This motivates the following
definition. 

\begin{Definition}[Dual of a definition of volume, \cite{alvarez_thompson}]
Let $\mu$ be a definition of volume. Then the dual definition of volume $\mu^*$ is given by 
\begin{displaymath}
 (\mu^*)_V:=(\mu_{V^*})^*.
\end{displaymath}
In terms of the associated functionals $\mathcal{V}$ and $\mathcal{V}^*$, we have 
\begin{displaymath}
 \mathcal{V}^*(B):=\frac{\vp(B)}{\mathcal{V}(B^\circ)}. 
\end{displaymath}
\end{Definition}

Examples
\begin{enumerate}
 \item Busemann's definition of volume and Holmes-Thompson's definition of volume are dual to each other.
\item Gromov's mass \eqref{item_gromov} and Gromov's mass* \eqref{item_benson} are dual to each other. 
\item Ivanov's definition of volume \eqref{item_ivanov} is dual to the definition of volume \eqref{item_dual_ivanov}.
\end{enumerate}

\begin{Definition}
\begin{enumerate}
 \item \label{item_k_density} Let $V$ be an $n$-dimensional vector space. A $k$-density on $V$, where $0 \leq k \leq n$,
is a map $\phi:\Lambda^k_s V \to \R$ such that $\phi(\lambda a)=|\lambda| \phi(a)$ for all $\lambda \in \R, a \in
\Lambda^k_s V$. Here $\Lambda^k_s V$ is the cone of simple $k$-vectors.
\item Let $M$ be an $n$-dimensional manifold. A $k$-density on $M$, where $0 \leq k \leq n$, is a continuous function
$\Phi:\Lambda^k_s(TM) \to \R$ such that the restriction to each tangent space $T_pM, p \in M$ is a $k$-density in the
sense of \eqref{item_k_density}. 
\end{enumerate}
\end{Definition}

\begin{Definition} 
Let $\mu$ be a definition of volume. If $V$ is a normed vector space, the induced $(n-1)$-volume
density $\mu_{n-1}$ on $V$ is defined as follows. If
$v_1 \wedge \ldots \wedge v_{n-1} \neq 0$, then 
\begin{displaymath} 
\mu_{n-1}(v_1 \wedge \ldots \wedge v_{n-1}):=\mu_W(v_1 \wedge \ldots \wedge v_{n-1}),
\end{displaymath}
where $W$ is the linear span of
these vectors with the induced norm. Otherwise $\mu_{n-1}(v_1 \wedge \ldots \wedge v_{n-1})=0$. 
\end{Definition}

\begin{Definition} \label{def_convex_volumedef}
A definition of volume $\mu$ is called
convex if $\mu_{n-1}:\Lambda^{n-1}V \to \R$ is a norm
for each normed space $V$.
\end{Definition}

Example: Busemann's definition of volume is convex. This is equivalent to the convexity of the intersection body (see
below). Holmes-Thompson volume is also convex. This is equivalent to the convexity of the projection body. Gromov's
mass* is convex, while Gromov's mass is not convex \cite{alvarez_thompson}. Ivanov's definition of volume is convex (\cite{ivanov08}, Thm. 6.2), while its dual is not convex. 

Given a compact hypersurface $M \subset V$, the norm on $V$ induces a Finsler metric on $M$ and we may integrate $\mu_{n-1}$ over $M$. We will be mostly interested in the case where $M=\partial K$ is the boundary of a smooth compact convex body. If $K$
is any compact convex body, then the integral still exists thanks to Alexandrov's theorem \cite{alexandrov39}. We will write
\begin{displaymath}
A_\mu(K):=\int_{\partial K} \mu_{n-1}
\end{displaymath}
and call $A_\mu(K)$ the surface area of $K$ with respect to $\mu$.

\begin{Proposition}[\cite{alvarez_thompson}] \label{prop_convexity}
 A definition of volume $\mu$ is convex if and only if for each pair $K \subset L$, we have 
\begin{displaymath}
 A_\mu(K) \leq A_\mu(L). 
\end{displaymath}
\end{Proposition}

\begin{Definition}[Isoperimetrix]
Let $\mu$ be a convex definition of volume. Then, for each normed space $(V,B)$ of dimension $n$ there exists a unique
centrally symmetric compact convex body $\mathbb{I}_\mu(B) \subset V$
such that for all compact convex bodies $K \subset V$ 
\begin{displaymath}
A_\mu(K) = n V(K[n-1],\mathbb{I}_\mu(B)).
\end{displaymath}
Here $V$ denotes the mixed volume \cite{schneider_book93}. The body $\mathbb{I}_\mu(B)$ is
called isoperimetrix. 
\end{Definition}

Strictly speaking, the isoperimetrix depends on the choice of a volume form. However, the definition of volume $\mu$
gives us a canonical choice of Lebesgue measure on $V$ which we will use in the following.  

Let us recall the construction of the isoperimetrix. The function $\mu_{n-1}:\Lambda^{n-1}V \to \R$ is convex and
$1$-homogeneous by assumption. The volume form on $V$ induces an isomorphism $\Lambda^{n-1}V \cong V^*$. We thus get a
convex and $1$-homogeneous function on $V^*$, which is the support function of the isoperimetrix. 

\begin{Proposition}[Isoperimetric inequality]
Let $B$ be the unit ball of a normed space and $\mu$ a convex definition of volume. Among all compact convex bodies $K$
with given volume, the surface
area with respect to $\mu$ is minimal precisely for a homothet of the isoperimetrix $\mathbb{I}_\mu(B)$.
Equivalently, for each choice of Lebesgue measure, we have 
\begin{displaymath}
 \frac{A_\mu(K)^n}{\mathcal{L}^n(K)^{n-1}} \geq \frac{A_\mu(\mathbb{I}_\mu(B))^n}{\mathcal{L}^n
(\mathbb{I}_\mu(B))^{n-1}}.
\end{displaymath}
\end{Proposition}

Examples: 
\begin{enumerate}
 \item The isoperimetrix for Busemann's volume is (up to a constant) the polar of the intersection body: 
\begin{displaymath}
 \mathbb{I}_{\mu^b}(B)=\omega_{n-1} I^\circ B.
\end{displaymath}

Let us briefly recall the definition of the intersection body, referring to \cite{haberl08, koldobsky, ludwig06, lutwak88, schuster08, thompson_book96} for details, more information and
generalizations. 

Given a non-zero volume form $\Omega \in \Lambda^n V^*$ and $\xi \in V^*, \xi \neq 0$, we may write (in a
non-unique way) $\Omega=\xi \wedge \Omega_\xi$ with $\Omega_\xi \in \Lambda^{n-1}V^*$. Then the restriction of
$\Omega_\xi$ to $\ker \xi$ is a volume form, which does not depend on the choice of $\Omega_\xi$. The intersection
body of a star body $S \subset V$ is the star body $IS$ in $V^*$ whose
radial function is given by 
\begin{displaymath}
 \rho(IS,\xi)=\vol(S \cap \ker \xi, \Omega_\xi).
\end{displaymath}
By a non-trivial result due to Busemann \cite{busemann49}, the intersection body of a centrally symmetric convex body is convex. 

Busemann and Petty \cite{busemann_petty56} have shown that the Busemann surface area $A_{\mu^b}(B)$ of the unit sphere is maximal
precisely if $B$ is a parallelotope. No lower bound seems to be known, except in dimensions $2$ and $3$. We refer to \cite{alvarez_thompson} and \cite{alvarez_thompson05} for more results in this direction. 
\item For the Holmes-Thompson definition of volume, the isoperimetrix is (up to a constant) the projection body of the
polar:
\begin{displaymath}
 \mathbb{I}_{\mu^{ht}}(B)=\frac{1}{\omega_{n-1}} \Pi(B^\circ).
\end{displaymath}

Let us also recall briefly the definition of the projection body. Again, we have to refer to the literature for a deeper
study of projection bodies \cite{koldobsky, ludwig02,schuster08, schuster_wannerer, thompson_book96}. 

If $v \in V, v \neq 0$, then $i_v\Omega:=\Omega(v,\cdot)$ is a
volume form on $V/\langle v\rangle$. Let $\pi_v:V \to V/\langle v\rangle$ be the projection map. The projection body
$\Pi K$ of a compact convex body $K$ is the compact convex body in $V^*$ whose support function is given by 
\begin{displaymath}
 h(\Pi K,v)=\vol(\pi_v K,i_v\Omega), \quad v \in V \setminus \{0\}. 
\end{displaymath} 
The projection body of a compact convex body is convex. 
\end{enumerate}

Let us recall a famous geometric inequality related to the projection body. 

\begin{Theorem}[Petty's projection inequality, \cite{petty71}] \label{thm_ppi}
Let $K \subset V$ be a compact convex body and $E \subset V$ an ellipsoid. Then 
\begin{displaymath}
 \vol(K)^{n-1} \vol \Pi^\circ K \leq \vol(E)^{n-1} \vol \Pi^\circ E 
\end{displaymath}
with equality precisely for ellipsoids. 
\end{Theorem}

Holmes and Thompson have shown that the Holmes-Thompson surface area $A_{\mu^{ht}}(B)$ is the same as the Holmes-Thompson
surface area $A_{\mu^{ht}}(B^\circ)$ (computed in the dual normed space $(V^*,B^\circ)$). This result was later
put in a symplectic geometry framework by Alvarez \cite{alvarez06}, who reproved their result and
showed that dual spheres have the same girth (the girth is the length of the shortest symmetric geodesic).
This confirmed a conjecture by Sch\"affer \cite{schaeffer_book}.

\section{The isoperimetrix in the dual Brunn-Minkowski theory}
\label{sec_dual_isoperimetrix}

The natural setting for the dual Brunn-Minkowski theory is that of star bodies (instead of convex bodies) and radial
addition (instead of Minkowski addition). We let $S_1 \tilde{+} S_2$ denote the radial addition of $S_1$
and $S_2$, i.e.
\begin{displaymath}
 \rho_{S_1 \tilde{+} S_2}=\rho_{S_1}+\rho_{S_2}, \quad S_1,S_2 \in \mathcal{S}(V),
\end{displaymath}
where $\rho$ denotes the radial function and $\mathcal{S}(V)$ is the space of star bodies. The natural topology on
$\mathcal{S}(V)$ is the radial topology. The dense subspace of star bodied with smooth radial function is denoted by
$\mathcal{S}^\infty(V)$.

Given an $n$-dimensional vector space $V$ with a volume form $\Omega \in \Lambda^n V^*$ and corresponding Lebesgue
measure $\mathcal{L}^n$, there is a unique functional                                               
\begin{displaymath}
 \tilde V:\mathcal{S}(V)^n \to \R
\end{displaymath}
called {\it dual mixed volume}, which is symmetric, multi-linear (with respect to
radial addition), continuous with respect to the radial topology and which satisfies $\tilde
V(S,\ldots,S)=\vol(S)$. It
was introduced by Lutwak \cite{lutwak75}.

Explicitly, we have for $S_1,\ldots,S_n \in \mathcal{S}(V)$ 
\begin{displaymath}
 \tilde V(S_1,\ldots,S_n)=\frac{1}{n} \int_{\Sigma} \rho_{S_1} \cdots \rho_{S_n} i_v \Omega,
\end{displaymath}
where $\Sigma$ is any $n-1$-submanifold in the same homology class as the sphere and $(i_v \Omega)|_v:=\Omega(v,\cdot)$
is an $(n-1)$-form on $V$. We refer to \cite{yang10} for more information on such {\it contour integrals}. In
particular, taking
$\Sigma:=\partial S$, 
\begin{displaymath}
 \tilde V(S[n-1],T)=\frac{1}{n} \int_{\partial S} \rho_T i_v \Omega.
\end{displaymath}

\begin{Theorem}[Dual Minkowski inequalities] \label{thm_dual_minkowski}
For $K,L \in \mathcal{S}(V)$ we have 
\begin{displaymath}
 \tilde V(K,\ldots,K,L)^n \leq \vol(K)^{n-1} \vol(L).
\end{displaymath}
\end{Theorem}

Let $S \in \mathcal{S}^\infty(V)$. Each tangent space $T_p \partial S$ may be considered as a subspace of $V$. But it
can also be considered as the quotient
\begin{displaymath}
T_p\partial S=V/\langle p\rangle.
\end{displaymath}
with the quotient norm. Explicitly, 
\begin{displaymath}
 \|v\|=\inf_{t \in \R} \|v+tp\|, \quad v \in T_p\partial S,
\end{displaymath}
which is not larger than the induced norm.

\begin{Definition} 
Let $\mu$ be a definition of volume. If $(V,B)$ is a normed vector space and $M \subset V$ an $(n-1)$-dimensional
submanifold such that $p \notin T_pM$ for all $p \in M$, there is an induced $(n-1)$-volume density $\tilde
\mu_{n-1}$ on $M$ defined as follows. 
Let $p \in M$ and $v_1,\ldots,v_{n-1} \in T_pM$. If $v_1 \wedge
\ldots \wedge v_{n-1} \neq 0$, then 
\begin{displaymath} 
\tilde \mu_{n-1}(v_1 \wedge \ldots \wedge v_{n-1}):=\mu_{T_pM}(v_1 \wedge \ldots \wedge v_{n-1}),
\end{displaymath}
where $T_pM=V/\langle p \rangle$ is endowed with the quotient norm. Otherwise $\mu_{n-1}(v_1 \wedge \ldots
\wedge v_{n-1})=0$. 
\end{Definition}

We may integrate $\tilde \mu_{n-1}$ over $M$. In the particular case where $M$ is the boundary of a star-shaped
set, this yields the following definition. 

\begin{Definition}
Let $\mu$ be a definition of volume and $S \in \mathcal{S}^\infty(V)$. We call 
\begin{displaymath}
 \tilde A_\mu(S):=\int_{\partial S} \tilde \mu_{n-1}
\end{displaymath}
the dual surface area of $S$ with respect to the definition of volume $\mu$.  
\end{Definition}

It will follow from Equation \eqref{eq_tilde_mu_as_mixed_volume} below that this
definition extends by continuity to all star bodies.

The quotient metric on the boundary of a unit ball was first studied by Faifman
\cite{faifman}, who showed an analogue of Alvarez' result mentioned in the last section. More precisely, he proved that
the girth, the length spectrum and Holmes-Thompson volume of the unit sphere and of the polar unit sphere are the same.

\proof[Proof of Theorem \ref{mainthm_dual_isop}]
In the following, the pairing between an element of a vector space and an element of the dual vector space will be
denoted by $\langle \cdot,\cdot\rangle$. 

Let $V$ be an $n$-dimensional normed space and $\Omega \in \Lambda^n V^*, \Omega \neq 0$ a volume form. There is an
isomorphism 
\begin{align*}
\tau: V & \to \Lambda^{n-1}V^*=(\Lambda^{n-1}V)^* \\
 w & \mapsto [a \mapsto \Omega(a \wedge w)].
\end{align*}
Let $v \in V, v \neq 0$. Clearly $\tau(v)$ vanishes on $(n-1)$-vectors of the form $a=v \wedge u, u \in \Lambda^{n-2}V$,
hence it belongs to $\Lambda^{n-1} \ker v$, where
$\ker v:=\{\eta \in V^*: \langle \eta,v\rangle=0\}$. Therefore $\langle v\rangle$ is mapped isomorphically to
$\Lambda^{n-1} \ker v$.

Let $W:=\left(V/\langle v\rangle,\pi_vB\right)$. We claim that there is a duality of normed spaces
\begin{displaymath}
W^*=\left(V/\langle v\rangle,\pi_vB\right)^*=\left(\ker v,B^\circ \cap \ker v\right).
\end{displaymath}
Indeed, the dual of the projection $\pi_v:V \twoheadrightarrow V/\langle v\rangle$ is the inclusion $\pi_v^*:\ker v
\hookrightarrow V^*$. 
Hence 
\begin{align*}
 (\pi_vB)^\circ & = \{\xi \in \ker v| \langle \xi,\pi_vp\rangle \leq 1 \forall p \in B\}\\
& =  \{\xi \in \ker v| \langle\pi_v^* \xi,p\rangle \leq 1 \forall p \in B\}\\
& = \{\xi \in \ker v| \pi_v^* \xi \in B^\circ\}\\
& = B^\circ \cap \ker v. 
\end{align*}

We define $h: V \to \R$ by $h(v):=\mu_{n-1}^*(\tau(v))$. Clearly $h$ is
$1$-homogeneous and positive. Since $\mu$ is the dual definition of volume of $\mu^*$, we have for $a \in \Lambda^{n-1}
W$ 
\begin{align*}
 \tilde \mu_{n-1}(a) & = \mu_W(a)\\
& = \sup\{\langle a,b\rangle|b \in \Lambda^{n-1}W^*, \mu_{W^*}^*(b) \leq 1\}\\
& =\sup\{\langle a, \tau(w)\rangle |w \in \langle v\rangle, \mu_{n-1}^*(\tau(w)) \leq 1\}\\
& =\sup\{\Omega(a \wedge w)|w \in \langle v\rangle, h(w) \leq 1\}\\
& = \frac{1}{h(v)} |\Omega(a \wedge v)|\\
& = \left|i_{\frac{v}{h(v)}} \Omega(a)\right|.
\end{align*}

Define the star body $T$ by 
\begin{equation} \label{eq_def_dual_iso}
 \rho_T(v)=\frac{1}{h(v)}. 
\end{equation}

It follows that 
\begin{equation} \label{eq_tilde_mu_as_mixed_volume}
 \tilde A_\mu(S) = \int_{\partial S} \tilde \mu_{n-1} = \int_{\partial S}  \rho_T(v) i_v \Omega
=  n \tilde V(S[n-1],T). 
\end{equation}
With $\tilde{\mathbb I}_\mu B:=T$, the first statement of the theorem follows. \\

By the dual Minkowski inequality (Theorem \ref{thm_dual_minkowski}), we find 
\begin{displaymath}
 \tilde A_\mu(S) \leq n \vol(S)^\frac{n-1}{n} \vol(\tilde{\mathbb I}_\mu B)^\frac{1}{n}
\end{displaymath}
with equality if and only if $S$ and $\tilde{\mathbb I}_\mu B$ are dilates of each other. This shows the second
statement of the theorem.\\

Finally, suppose that $\mu^*$ is convex in the sense of Definition \ref{def_convex_volumedef}. By the construction of
the isoperimetrix, the function $h$ defined above is the support
function of $\mathbb{I}_{\mu^*} B^\circ$. From \eqref{eq_relation_support_radial} and \eqref{eq_def_dual_iso} it follows
that $\tilde{\mathbb{I}}_{\mu} B$ is
the polar of $\mathbb{I}_{\mu^*} B^\circ$, which finishes the proof.  
\endproof

\proof[Proof of Theorem \ref{mainthm_surface_sphere}]
Let $E$ be an ellipsoid. By the dual Minkowski inequality (Theorem \ref{thm_dual_minkowski}) and Petty's projection
inequality
(Theorem \ref{thm_ppi}) we have 
\begin{align*}
\tilde A_{{\mu}^b}(B) & = n \omega_{n-1} \tilde V(B[n-1],\Pi^\circ B) \\
& \leq  n \omega_{n-1}  (\vol B)^\frac{n-1}{n}
(\vol \Pi^\circ B)^\frac{1}{n} \\
& \leq  n \omega_{n-1}  (\vol E)^\frac{n-1}{n}
(\vol \Pi^ \circ E)^\frac{1}{n} \\
& = n \omega_{n-1} \tilde V(E[n-1],\Pi^ \circ E) \\
& =\tilde A_{{\mu}}(E)\\
& = n \omega_n.
\end{align*}

The equality case follows from the fact that equality in Petty's projection inequality is
attained for (not necessarily centered)
ellipsoids and that $B=-B$.
\endproof

\proof[Proof of Corollary \ref{cor_max_quot_girth}]
Clearly the quotient girth of $B$ is not larger than the quotient girth of any central two-dimensional section, which is
bounded by $2\pi$ by the above theorem. If $B$ is not an ellipsoid, then there exists a two-dimensional central section
which is not an ellipsoid (see \cite{gardner_book06}, Thm. 7.1.5) and whose quotient girth is strictly smaller than $2\pi$.  
\endproof


\begin{thebibliography}{10}

\bibitem{alexandrov39}
A.~D. Alexandroff.
\newblock Almost everywhere existence of the second differential of a convex
  function and some properties of convex surfaces connected with it.
\newblock {\em Leningrad State Univ. Annals [Uchenye Zapiski] Math. Ser.},
  6:3--35, 1939.

\bibitem{alvarez06}
J.~C. {\'A}lvarez~Paiva.
\newblock Dual spheres have the same girth.
\newblock {\em Amer. J. Math.}, 128(2):361--371, 2006.

\bibitem{alvarez_thompson}
J.~C. {\'A}lvarez~Paiva and A.~C. Thompson.
\newblock Volumes on normed and {F}insler spaces.
\newblock In {\em A sampler of {R}iemann-{F}insler geometry}, volume~50 of {\em
  Math. Sci. Res. Inst. Publ.}, pages 1--48. Cambridge Univ. Press, Cambridge,
  2004.

\bibitem{alvarez_thompson05}
Juan~Carlos {\'A}lvarez~Paiva and Anthony Thompson.
\newblock On the perimeter and area of the unit disc.
\newblock {\em Amer. Math. Monthly}, 112(2):141--154, 2005.

\bibitem{ambrosio_kirchheim00}
Luigi Ambrosio and Bernd Kirchheim.
\newblock Currents in metric spaces.
\newblock {\em Acta Math.}, 185(1):1--80, 2000.

\bibitem{busemann_petty56}
H.~Busemann and C.~M. Petty.
\newblock Problems on convex bodies.
\newblock {\em Math. Scand.}, 4:88--94, 1956.

\bibitem{busemann49b}
Herbert Busemann.
\newblock The isoperimetric problem for {M}inkowski area.
\newblock {\em Amer. J. Math.}, 71:743--762, 1949.

\bibitem{busemann49}
Herbert Busemann.
\newblock A theorem on convex bodies of the {B}runn-{M}inkowski type.
\newblock {\em Proc. Nat. Acad. Sci. U. S. A.}, 35:27--31, 1949.

\bibitem{faifman}
Dmitry Faifman.
\newblock An extension of {S}ch\"affer's dual girth conjecture to
  {G}rassmannians.
\newblock {\em J. Differential Geom.}, 92(1):201--220, 2012.

\bibitem{gardner_book06}
Richard~J. Gardner.
\newblock {\em Geometric tomography}, volume~58 of {\em Encyclopedia of
  Mathematics and its Applications}.
\newblock Cambridge University Press, Cambridge, second edition, 2006.

\bibitem{haberl08}
Christoph Haberl.
\newblock {$L_p$} intersection bodies.
\newblock {\em Adv. Math.}, 217(6):2599--2624, 2008.

\bibitem{ivanov08}
S.~V. Ivanov.
\newblock Volumes and areas of {L}ipschitz metrics.
\newblock {\em Algebra i Analiz}, 20(3):74--111, 2008.

\bibitem{koldobsky}
Alexander Koldobsky.
\newblock {\em Fourier analysis in convex geometry}, volume 116 of {\em
  Mathematical Surveys and Monographs}.
\newblock American Mathematical Society, Providence, RI, 2005.

\bibitem{ludwig02}
Monika Ludwig.
\newblock Projection bodies and valuations.
\newblock {\em Adv. Math.}, 172(2):158--168, 2002.

\bibitem{ludwig06}
Monika Ludwig.
\newblock Intersection bodies and valuations.
\newblock {\em Amer. J. Math.}, 128(6):1409--1428, 2006.

\bibitem{ludwig10}
Monika Ludwig.
\newblock Minkowski areas and valuations.
\newblock {\em J. Differential Geom.}, 86(1):133--161, 2010.

\bibitem{lutwak75}
Erwin Lutwak.
\newblock Dual mixed volumes.
\newblock {\em Pacific J. Math.}, 58(2):531--538, 1975.

\bibitem{lutwak88}
Erwin Lutwak.
\newblock Intersection bodies and dual mixed volumes.
\newblock {\em Adv. in Math.}, 71(2):232--261, 1988.

\bibitem{lutwak_handbook93}
Erwin Lutwak.
\newblock Selected affine isoperimetric inequalities.
\newblock In {\em Handbook of convex geometry, {V}ol.\ {A}, {B}}, pages
  151--176. North-Holland, Amsterdam, 1993.

\bibitem{morgan92}
Frank Morgan.
\newblock Minimal surfaces, crystals, shortest networks, and undergraduate
  research.
\newblock {\em Math. Intelligencer}, 14(3):37--44, 1992.

\bibitem{petty71}
Clinton~M. Petty.
\newblock Isoperimetric problems.
\newblock In {\em Proceedings of the {C}onference on {C}onvexity and
  {C}ombinatorial {G}eometry ({U}niv. {O}klahoma, {N}orman, {O}kla., 1971)},
  pages 26--41. Dept. Math., Univ. Oklahoma, Norman, Okla., 1971.

\bibitem{schaeffer_book}
Juan~Jorge Sch{\"a}ffer.
\newblock {\em Geometry of spheres in normed spaces}.
\newblock Marcel Dekker Inc., New York, 1976.
\newblock Lecture Notes in Pure and Applied Mathematics, No. 20.

\bibitem{schneider_book93}
Rolf Schneider.
\newblock {\em Convex bodies: the {B}runn-{M}inkowski theory}, volume~44 of
  {\em Encyclopedia of Mathematics and its Applications}.
\newblock Cambridge University Press, Cambridge, 1993.

\bibitem{schuster08}
Franz~E. Schuster.
\newblock Valuations and {B}usemann-{P}etty type problems.
\newblock {\em Adv. Math.}, 219(1):344--368, 2008.

\bibitem{schuster_wannerer}
Franz~E. Schuster and Thomas Wannerer.
\newblock {${\rm GL}(n)$} contravariant {M}inkowski valuations.
\newblock {\em Trans. Amer. Math. Soc.}, 364(2):815--826, 2012.

\bibitem{taylor73}
Jean~E. Taylor.
\newblock Unique structure of solutions to a class of nonelliptic variational
  problems.
\newblock In {\em Differential geometry ({P}roc. {S}ympos. {P}ure. {M}ath.,
  {V}ol. {XXVII}, {S}tanford {U}niv., {S}tanford, {C}alif., 1973), {P}art 1},
  pages 419--427. Amer. Math. Soc., Providence, R.I., 1975.

\bibitem{taylor78}
Jean~E. Taylor.
\newblock Crystalline variational problems.
\newblock {\em Bull. Amer. Math. Soc.}, 84(4):568--588, 1978.

\bibitem{thompson_book96}
Anthony~C. Thompson.
\newblock {\em Minkowski geometry}, volume~63 of {\em Encyclopedia of
  Mathematics and its Applications}.
\newblock Cambridge University Press, Cambridge, 1996.

\bibitem{wulff01}
G.~Wulff.
\newblock Zur {F}rage der {G}eschwindigkeit des {W}achstums und der
  {A}ufl\"osung der {K}rystallfl\"achen.
\newblock {\em Z. Krystallogr. Mineral.}, 34:449--530, 1901.

\bibitem{yang10}
Deane Yang.
\newblock Affine geometry using the homogeneous contour integral.
\newblock In {\em Proceedings of the 14th {I}nternational {W}orkshop on
  {D}ifferential {G}eometry and the 3rd {KNUGRG}-{OCAMI} {D}ifferential
  {G}eometry {W}orkshop [{V}olume 14]}, pages 1--12. Natl. Inst. Math. Sci.
  (NIMS), Taej\u on, 2010.

\end{thebibliography}
\end{document}